\input amstex 
\documentstyle{amsppt}
\nologo
\NoBlackBoxes

\magnification 1200
\hsize= 14.6 truecm
\vsize= 19.25 truecm
\hoffset= .85 truecm
\voffset= 2.1 truecm

\def\ns{\hskip -.15 truecm}

\def\preone{1.1 }
\def\pretwo{1.2 }
\def\prethree{1.3 }
\def\prefour{1.4 }
\def\main{Theorem 2.1 }
\def\maincurve{Proposition 2.2 }
\def\split{Lemma 2.3 }
\def\maincurvepr{2.4 }
\def\algstr{2.4.1 }
\def\maincurvecor{Corollary 2.6 }
\def\mainpr{2.7 }
\def\maincor{Corollary 2.8 }

\def\nocyc{Proposition 3.1 }
\def\nocyccor{Corollary 3.2 }
\def\noodd{Proposition 3.3 }
\def\Fzero{Example 3.4 }
\def\Fone{Example 3.5 }
\def\cone{Example 3.7 }

\def\hyp{Theorem 4.1 }
\def\Noet{Theorem 4.2 }
\def\Noetcor{Corollary 4.3 }

\def\CYth{Theorem 5.1 }

\def\Be{[Be] }
\def\Bo{[Bo] }
\def\Ca{[Ca] }
\def\Ci{[Ci] }
\def\EH{[EH] }
\def\F{[F] }
\def\sos{[GP1] }
\def\CY{[GP2] }
\def\can{[GP3] }
\def\trig{[GP4] }
\def\G{[G] }
\def\Hoone{[H1] }
\def\Hotwo{[H2] }
\def\Hothree{[H3] }
\def\Hofour{[H4] }
\def\HM{[HM] }
\def\Kod{[Kod] }
\def\Kon{[Kon] }
\def\M{[M] }
\def\OP{[OP] }
\def\T{[T] }

\topmatter
\title On the canonical ring of covers of surfaces of minimal degree
\endtitle
\author Francisco Javier Gallego \\ and \\ B. P. Purnaprajna
\endauthor
\address{Francisco Javier Gallego: Dpto. de \'Algebra,
 Facultad de Matem\'aticas,
 Universidad Complutense de Madrid, 28040 Madrid,
Spain}\endaddress
\email{gallego\@eucmos.sim.ucm.es}\endemail
\address{ B.P.Purnaprajna:
405 Snow Hall,
  Dept. of Mathematics,
  University of Kansas,
  Lawrence, Kansas 66045-2142
}\endaddress
\email{purna\@math.ukans.edu}\endemail
\endtopmatter
\document 

\vskip .3 cm

\headline={\ifodd\pageno\rightheadline \else\leftheadline\fi}
\def\rightheadline{\tenrm\hfil \eightpoint ON THE CANONICAL RING OF
  COVERS OF 
  SURFACES OF MINIMAL DEGREE
 \hfil\folio}
\def\leftheadline{\tenrm\folio\hfil \eightpoint F.J. GALLEGO \&
  B.P. PURNAPRAJNA \hfil}

\heading  Introduction \endheading

Let $\varphi :X \rightarrow Y$ be a generically finite morphism. The
purpose of this paper is to show how the $\Cal O_Y$-algebra 
structure on $\varphi_{*}\Cal{O}_{X}$ controls algebro-geometric aspects
of $X$ like the ring generation of graded rings associated to $X$ 
and the very ampleness of line bundles on $X$. 
As the main application of this we prove some new results for certain
regular surfaces $X$ of general type. Precisely, we find the degrees
of the generators of the canonical ring of $X$ when the canonical
morphism of $X$ is a finite cover of a surface of minimal degree. 
These results 
complement results of Ciliberto \Ci and Green \G \ns. The techniques
of this paper also yield 
different proofs of some earlier results, such as Noether's theorem
for certain kind of curves and some results on Calabi-Yau threefolds
that had appeared in \CY \ns.

\medskip

The canonical ring of surfaces of general type have
attracted the attention of several geometers. Kodaira (see \Kod \ns)
first proved that 
$|K_{X}^{\otimes m}|$ embeds a minimal surface of general type $X$ as 
a projectively normal variety for all $m\geq 8$. This was later
improved by 
Bombieri (see \Bo \ns) who proved the same
result if $m\geq 6$ and by Ciliberto (see \Ci \ns), who lowered the bound to
$m\geq
5$.
Recently the
authors proved (see \sos \ns) more general results on projective normality
and
higher syzygies for adjunction bundles for an algebraic surface. As a
corollary of these results they recovered and improved the results of
Bombieri and Ciliberto on projective normality and extended them to
higher syzygies.

An important class of minimal surfaces of general type are those
whose canonical divisor is base-point-free. 
If one
goes by the results proved so far surfaces with base-point-free
canonical bundle fall into one of these two categories: 
those whose canonical morphism maps onto a
surface of minimal degree and those whose canonical morphism doesn't
map onto a
surface of minimal degree. 
The latter  have been studied by Ciliberto (see \Ci \ns)
and (see Green \G \ns). The former are studied in this article.  Green
and Ciliberto proved 
this nice result regarding the generators of the canonical ring of $X$:

\medskip

Let $X$ be a regular surface of general type with a base point free
canonical divisor.
Assume that the canonical morphism $\varphi $ satisfies the following
conditions:

(1)$\varphi $ does not map $X$ generically $2:1$ onto $\bold{P}^{2}.$

(2) $\varphi (S)$ is not a surface of minimal degree (other than
$\bold{P}^{2}).$

Then the the canonical ring of $X$ is generated in degree less than or
equal to $2$.

\medskip 

In the present article we deal with surfaces 
of general type $X$ 
 whose canonical morphism $\varphi$ maps
$X$ onto a surface of minimal degree $Y$. 
These surfaces have been studied in the works of Horikawa (see \Hoone,
\Hotwo, \Hothree and \Hofour \ns), Catanese ((see \Ca \ns) and Konno
(see \Kon \ns) among others, 
where they play a
central role in the classification of surfaces of general type with small
$c_1^2$, in 
questions about degenarations and the moduli of surfaces of general
part.
The study of these surfaces have a direct bearing on the study of linear
series on threefolds such as Calabi-Yau threefolds as the results in
\OP  and
authors results in \CY  show.

The study the canonical rings of these surfaces is carried out in Section
2.
We determine the precise degrees of the generators of its  canonical
ring (see \main \ns). The answer depends on the degree of
$\varphi$ and the degree of $Y$.  As a corollary
of our result and the result of Ciliberto and Green,
we find that conditions (1) and (2) above characterize the regular
surfaces of general type with base
point free canonical bundle whose canonical ring is generated in
degree
less than or equal to $2$.

\medskip

As we said at the beginning, in order to study the canonical ring of 
$X$ we will use the existence of a generically finite morphism
$\varphi$ from $X$ to 
a variety $Y$.   The morphism $\varphi$ is the canonical morphism of
$X$ and $Y$ is a surface of minimal degree, that
is,
a nondegenerate surface in projective space whose degree is equal to its
codimension plus $1$. The classification of these surfaces is classically
known: they are (linear) $\bold P^{2}$, the Veronese surface in $\bold
P^{5}$, smooth rational scrolls or cones over one of them (see
\EH \ns). Thus the surface $Y$ is simpler than $X$. A
measurement of its simplicity is that its general hyperplane section is a
smooth, rational normal curve. 
Therefore in
Section 2 we see how the algebra structure of $\varphi _{*}\Cal{O}_{X}$
governs the multiplicative structure of the canonical ring of $X$ and
we use this to study its ring generators. 

\medskip

Section 3 is devoted to constructing examples. We recall some known
examples of surfaces of general type mapping to a surface of minimal
degree and construct some new ones. We also show that certain kinds of
examples of finite canonical morphisms are not possible. For example, we
show that odd degree covers of smooth rational scrolls or cyclic covers of
degree bigger than 3 of surfaces of minimal degree do not exist.

\medskip

In Section 4, we show another example that illustrates the relation
between the algebra structure of $\pi _{*}\Cal{O}_{X}$ given by a
finite morphism $\pi$ and the canonical
ring of $X$. As an application of this, we
give a different proof
of Noether's theorem for curves general among those possessing an
effective theta-characteristic. In this case the finite morphism $\pi$
is
induced
by the complete linear series of a theta-characteristic of a curve mapping
to a rational curve. In Section 5, we give yet another illustration of
this philosophy  and give a different proof of results on Calabi-Yau
threefolds proved in \CY \ns. 

\medskip

Finally we expand on these ideas in two forthcoming articles, \can and
\trig \ns. In the first we study the canonical ring of higher dimensional
varieties of general type whose canonical morphism maps onto a variety of
minimal degree. One of results in \can shows that the converse of the
theorem of Ciliberto and Green for surfaces stated above is false for
higher
dimensional varieties of general type. In the second we carry out a
detailed
study of homogeneous rings associated to line bundles on trigonal curves.

\proclaim{Convention} We will work over an algebraic closed field of
characteristic $0$ \endproclaim

\heading 1. Preliminaries \endheading

In this section we will recall some known facts about the push forward
of the structure sheaf of a variety by a flat, 
finite morphism. We summarize these facts below and refer for the proof to
\HM \ns, Section 2.  

Let $X$ and $Y$ be algebraic varieties over
a field $\bold k$ and let $n$ a natural number which does not divide
car$(\bold k)$. 
Let  $\pi: X
\longrightarrow Y$ be a finite, flat morphism of degree $n$. We have
the folowing facts:

\medskip

\item{\bf \preone\ns.} The sheaf
$\pi_* \Cal O_X$ is a rank $n$ locally free sheaf on $Y$  of algebras
over $\Cal O_Y$.

\medskip  
 
\item{\bf \pretwo\ns.} There exists a map  
$$\frac 1 n \text{tr}: \pi_* \Cal O_X \longrightarrow \Cal O_Y$$
of sheaves of $\Cal
  O_Y$ modules defined locally as follows:
Given $\alpha \in \pi_*\Cal O_X$ we consider the 
homomorphism of $\Cal O_Y$-modules 
$$ \pi_*\Cal O_X @>\cdot  \alpha >> \pi_*\Cal O_X $$
induced by multiplication by $\alpha$.
Then we define $\frac 1 n \text{tr} (\alpha)$ as the trace
of such homomorphism divided by $n$.

\medskip

\item{\bf \prethree\ns.} $\frac 1 n \text{tr}$ is surjective, in
  fact, the map  
  $\Cal O_Y \hookrightarrow \pi_*\Cal O_X$ induced by $\pi$ 
is a section of $\frac 1 n \text{tr}$. Therefore the sequence
$$0 \longrightarrow E \longrightarrow \pi_* \Cal O_X @>\frac 1 n
\text{tr} >>  \Cal O_Y$$ 
splits.   
$E$ is the kernel of $ \frac 1 n \text{tr}$ and locally  consists of
the trace $0$ elements of $\pi_*\Cal O_X$. We will call $E$ the {\it
  trace-zero module} of $\pi$. 

\medskip

\item{\bf \prefour\ns.} $\pi_*(\Cal O_X)$ is a sheaf of $\Cal
  O_Y$-algebras, therefore 
  it has a multiplicative structure. Its multiplication map is an
  $\Cal O_Y$-bilinear map
$$[\Cal O_Y \oplus E] \otimes [\Cal O_Y \oplus E] \longrightarrow
\Cal O_Y \oplus E $$
made of four components. 
The first component 
$$\Cal O_Y \otimes \Cal O_Y \longrightarrow \Cal O_Y \oplus E$$ 
is given by the multiplication in $\Cal O_Y$ and therefore goes to
$\Cal O_Y$. 
The components
$$\displaylines{\Cal O_Y \otimes E \longrightarrow \Cal O_Y \oplus E
  \cr
E \otimes \Cal O_Y  \longrightarrow \Cal O_Y \oplus E}$$
are given by the left and right module structure of $E$ over $\Cal
O_Y$ and therefore go to $E$. 
Finally there is a fourth component
$$E \otimes E \longrightarrow \Cal O_Y \oplus E$$
which factorsizes through
$$S^2E \longrightarrow \Cal O_Y \oplus E$$
for multiplication in $\pi_*\Cal O_Y$ is commutative.

\heading 2. Covers of surfaces of minimal
degree. \endheading

Our purpose is to study the generators of the canonical ring of
certain surfaces of general type. Precisely we are interested
in studying 
those regular surfaces of general type whose 
canonical divisor is  base point free and such that the image of the
canonical morphism is a variety of minimal degree. We obtain the
following

\proclaim{\main} Let $S$ be a regular surface of general
type with at worst canonical singularities and such that its canonical
bundle $K_S$ is base-point-free. Let
$\varphi$ be the canonical morphism of $S$. Let $n$ be the
degree of $\varphi$ and assume that the image of $\varphi$ is a
surface of minimal degree $r$. Then

\item{1)} if $n=2$ and $r=1$ (i.e., if $\varphi$ is generically
$2:1$ onto $\bold P^2$), the canonical ring of $S$ is generated by
its part of degree $1$ and one generator in degree $4$;

\item{2)} if $n \neq 2$ or $r \neq 1$, the canonical ring of $S$ is
generated by its part of degree $1$, 
$r(n-2)$ generators in degree $2$ and $r-1$ 
generators in degree $3$.  

\endproclaim

The knowledge of how many linearly independent generators are needed
in each degree is obtained from the knowledge of the image of the
multiplication maps of global sections of powers of the canonical
bundle. We study those multiplication maps by studying similar maps of
a curve $C$ in $|K_S| $. Thus we will first prove the following 

\proclaim{\maincurve} Let $C$ be a smooth curve.  Let $\theta$
be a base-point-free line bundle on $C$ such that
$\theta^{\otimes 2}=K_C$. Let
$\pi$ be the morphism induced by $|\theta|$, let $n$ be the
degree of
$\pi$ and assume that
$\pi(C)$ is a rational normal curve of degree $r$. 
Let $\beta(s,t)$ be the multiplication map 
$$ H^0(\theta^{\otimes s}) \otimes H^0(\theta^{\otimes t})
\longrightarrow  H^0(\theta^{\otimes s+t})  \text{, for all }
s, t >0 \ .$$ The codimension of the image of $\beta(s,t)$ in
$H^0(\theta^{\otimes s+t})$ is as follows:

\item{a)} If $r=1$, the codimension is:
\itemitem{a.1)} $n-2$, for $s=t=1$,
\itemitem{a.2)} $0$, for $s=2, t=1$, i.e., $\beta(2,1)$ surjects.
\itemitem{a.3)} $1$, for $s=3, t=1$.
\itemitem{a.4)} $1$, for $s=t=2$, $n=2$ and $0$ if $n >2$.  
\itemitem{a.5)} $0$, for $s \geq 4, t=1$, i.e., $\beta(s,1)$ surjects for
all $s \geq 4$.

\item{b)} If $r > 1$, the codimension is: 
\itemitem{b.1)} $r(n-2)$, for $s=t=1$. 
\itemitem{b.2)} $r-1$, for $s=2, t=1$,
\itemitem{b.3)} $0$, for $s \geq 3, t=1$, i.e., $\beta(s,1)$ surjects for
all $s \geq 3$.

Moreover, if $r=1$ and $n=2$, then the image of $\beta(2,2)$ and the image of
$\beta(3,1)$ are equal. 

\endproclaim

In order to prove \maincurve we will use the following

\proclaim{\split} Let $C$, $\theta$ and $\pi$ as in the statement of
\maincurve \ns. Then
$$\pi_*\Cal O_C=\Cal O_{\bold P^1} \oplus (n-2) \Cal O_{\bold
P^1}(-r-1) \oplus \Cal O_{\bold 
P^1}(-2r-2) \ .
$$
\endproclaim

{\it Proof.}
Since the image of $\pi$ is smooth and of
dimension $1$, $\pi$ is flat. Then 
$\pi_*\Cal O_C=\Cal O_{\bold P^1}
\oplus E$ as $\Cal O_{\bold P^1}$-modules, with $E$ vector bundle
over $\bold P^1$.  We now show that
$$E= (n-2)\Cal O_{\bold P^1}(-r-1) \oplus \Cal O_{\bold
P^1}(-2r-2) \ .
$$ 
We have $\pi_* \theta=\pi_*\Cal O_C\otimes \Cal O_{\bold
P^1}(r)$ and
$\pi_* K_C = \pi_*\Cal O_C\otimes \Cal O_{\bold P^1}(2r)$, by
projection formula. Any vector bundle over $\bold P^1$
splits, hence $$\pi_*(\Cal O_C)= \Cal O_{\bold P^1} \oplus E=
\Cal O_{\bold P^1} 
\oplus
\Cal O_{\bold P^1}(a_1) \oplus \cdots \oplus \Cal O_{\bold
P^1}(a_{n-1}) \ ,$$ for some negative integers $a_1, \dots,
a_{n-1}$ ($C$ is connected). Then
$h^1(K_C)=1$ implies that exactly one of the
$a_i$s, let us say
$a_{n-1}$, satisfies $a_{n-1}+2r=-2$. On the other hand, since 
$\pi$ is induced by the complete linear series $|\theta|$,
$h^0(\theta)=r+1= h^0(\Cal O_{\bold P^1}(r))$, so $a_i+r \leq -1$
for all $1 \leq i \leq n-2$. Finally, since degree of $\theta$ is
$g(C)-1$,
$h^1(\theta)=h^0(\theta)=r+1$. Since $h^1(\Cal O_{\bold
P^1}(-r-2))=r+1$, $a_i+r \geq -1$ for all $1 \leq i \leq n-2$,
so 
$a_i+r = -1$ for all $1 \leq i \leq n-2$. $\square$

\medskip

(\maincurvepr \ns) {\it Proof of \maincurve \ns.}
In \split we have completely determined the structure of $\pi_*\Cal
O_C$ as $\Cal O_{\bold P^1}$-module. Now we look at the structure
of $\pi_*\Cal O_C$ as $\Cal O_{\bold P^1}$-algebra. If $n=2$, it is 
completely determined by the branch divisor of $\pi$ on $\bold
P^1$, since in this case $\pi$ is cyclic. 
If $n >2$, we observe the following:
$$\displaylines{ \text{ (\algstr \ns) For
some $1 \leq i,j \leq n-2$, the projection of the map } \cr
\Cal
O_{\bold P^1}(a_i) \otimes \Cal O_{\bold P^1}(a_j) \longrightarrow
\pi_*\Cal O_C \text{ to }\Cal O_{\bold P^1}(-2r-2) \cr
\text{ is
surjective, in fact, it is an isomorphism.}}$$   

This is so because otherwise $\Cal O_{\bold P^1}
\oplus \Cal O_{\bold P^1}(a_1)
\oplus \cdots \oplus \Cal O_{\bold P^1}(a_{n-2})$ would be an
integral subalgebra of $\pi_*\Cal O_C$, free over $\Cal
O_{\bold P^1}$ of rank $n-1$. Then
$n-1$ should divide $n$, which is not possible if $n >2$.

\medskip

Now we will use our knowledge of $\pi_*\Cal O_X$ to study the maps
$\beta(s,r)$ which appear in the statement of the  proposition. We
will write $\beta_s$ in place of $\beta(s,1)$.  Let 
$R_l=H^0(\theta^{\otimes l})$. 
Then, since $\theta=\pi^*\Cal O_{\bold P^1}(r)$,  by projection formula
$$\displaylines{R_1 = H^0(\Cal O_{\bold P^1}(r)), \cr 
R_l= H^0(\Cal
O_{\bold P^1}(lr)
\oplus (n-2)H^0(\Cal O_{\bold P^1}((l-1)r-1)) \oplus H^0(\Cal
O_{\bold P^1}((l-2)r-2)) \text{ and }
\cr R_{l+1}= H^0(\Cal O_{\bold P^1}((l+1)r))
\oplus (n-2)H^0(\Cal O_{\bold P^1}(lr-1)) \oplus H^0(\Cal O_{\bold
P^1}((l-1)r-2))
\ . }$$
Therefore an element of $R_l$, i.e., a global section of
$H^0(\theta^{\otimes l})$ is a sum of $n$ components,
one in each piece of the above decomposition of $R_l$. On the
other hand, the product of an element of $R_l$ belonging to one of
the blocks with an element of $R_1$ is determined by the ring
structure of $\Cal O_{\bold P^1}$ and by the module structure of
$E$. More precisely, the restriction of $\beta_l$ to $H^0(\Cal
O_{\bold P^1}(lr)) \otimes H^0(\Cal O_{\bold P^1}(r))$ maps, in
fact isomorphically, onto $H^0(\Cal O_{\bold P^1}((l+1)r)$. The
restriction of $\beta_l$ to each of the blocks $H^0(\Cal
O_{\bold P^1}((l-1)r-1)) \otimes H^0(\Cal O_{\bold P^1}(r))$ maps
to the corresponding $H^0(\Cal O_{\bold P^1}(lr-1))$. This
restriction is $0$ if $(l-1)r-1$ is negative and an isomorphism
otherwise. Likewise, the restriction of $\beta_l$ to $H^0(\Cal
O_{\bold P^1}((l-2)r-2)) \otimes H^0(\Cal O_{\bold P^1}(r))$ goes
to $H^0(\Cal O_{\bold
P^1}((l-1)r-2))$, being $0$ if $(l-2)r-2$ is negative and an
isomorphism otherwise. Therefore it is crucial to tell which
blocks of a given $R_l$ are $0$. We have
$$\displaylines{
R_1=H^0(\Cal
O_{\bold P^1}(r)) ; \cr 
R_2=H^0(\Cal O_{\bold P^1}(2r)) \oplus
(n-2)H^0(\Cal O_{\bold P^1}(r-1)); \text{ and }
\text{ if } l \geq 3, \cr
 R_l=H^0(\Cal O_{\bold P^1}(lr)) \oplus  (n-2)H^0(\Cal O_{\bold
P^1}((l-1)r-1)) \oplus H^0(\Cal O_{\bold
P^1}((l-2)r-2))\ .}$$ 
All the direct summands appearing in the above
formulae are nonzero, except \linebreak 
$H^0(\Cal O_{\bold
P^1}((l-2)r-2))$ when $l=3$ and $r=1$ and $(n-2)H^0(\Cal O_{\bold
P^1}((l-1)r-1))$ for all $l$ and all $r$ when $n=2$.
We now determine the image of $\beta_l$. If $l=1$, the
image  of
$\beta_1$ is
$H^0(\Cal O_{\bold P^1}(2r))$, which has codimension $(n-2)r$ in
$R_2$. If $l=2$, the image  of $\beta_2$ is $H^0(\Cal
O_{\bold P^1}(3r)) \oplus H^0( (n-2)\Cal O_{\bold P^1}(2r-1))$ which has
codimension
$r-1$ in $R_3$. If $l =3$ and $r
\geq 2$ or if $l \geq 4$, the image  of $\beta_l$ is
all $R_l$, i.e., $\beta_l$ surjects. 
All this proves a.1), a.2), a.5) and b). If $r=1$  the image
of $\beta(3,1)$ is $H^0(\Cal O_{\bold P^1}(4r)) \oplus (n-2)H^0(\Cal
O_{\bold P^1}(3r-1))$, which has codimension $1$ in $R_4$. This proves
a.3). If $r=1$ and $n=2$, the image of $\beta(2,2)$  is $H^0(\Cal
O_{\bold P^1}(4r))$, which has codimension
$1$ in $R_4$. This proves the first claim in a.4) and the last
sentence of \maincurve \ns. Finally, if
$n >2$, recall (see \algstr \ns) that for some $1 \leq i,j \leq n-2$, the
projection of the map $$\Cal O_{\bold P^1}(a_i) \otimes \Cal
O_{\bold P^1}(a_j) \longrightarrow \pi_*\Cal O_C $$ to $\Cal
O_{\bold P^1}(-4)$ is surjective, in fact, it is an isomorphism.  
Then if $n >2$ the image of $\beta(2,2)$ is all $R_4$. This proves the
second part of a.4).
$\square$

\medskip

\noindent{\bf Remark 2.5.} Note that $\theta^{\otimes 2}=K_C$. Then a
proof of a.4), alternate to the one given above, 
can be obtained from Noether's  Theorem and 
from the base-point-free pencil
trick. 
The way how Noether's theorem is related to the algebra structure of
$\pi_*\Cal O_C$ will be clear in Section 4, where we will give
a different, simple proof of this classical result in certain particular
cases.

\medskip  

>From \maincurve we obtain the following

\proclaim{\maincurvecor} Let $C$ be a smooth curve. Let $\theta$
be a base-point-free line bundle on $C$ such that
$\theta^{\otimes 2}=K_C$. Let
$\pi$ be the morphism induced by $|\theta|$, let $n$ be the
degree of
$\pi$ and assume that
$\pi(C)$ is a rational normal curve. Let $R$ be 
$\bigoplus_{l=0}^\infty H^0(\theta^{\otimes l})$. Then 

\item{1)} if $r=1$ and $n=2$, the ring $R$ is
generated by its part of degree $1$ and one  generator in
degree $4$;

\item{2)} if  $r = 1$ and $n > 2$,  the ring $R$ is
generated by its part of degree $1$ and $n-2$  generators in
degree $2$;

\item{3)} if  $r > 1$,  the ring $R$ is
generated by its part of degree $1$, $r(n-2)$ generators in
degree $2$ and $r-1$ generators in degree $3$.
\endproclaim

{\it Proof:} To know in what degrees we need  generators we look
at the maps $\beta(s,t)$ of multiplication of sections. Precisely 
the number of generators
needed in degree $l+1$ is the codimension in $R_{l+1}$ of the sum of
the images of $\beta(l,1),\beta(l-1,2),\dots, \beta(\lfloor
  \frac{l+1}{2} \rfloor, \lceil \frac{l+1}{2} \rceil)$. In particular
$R$ is generated in degree less than or equal to $l$ if $\beta_{k}$
surjects for all $k \geq l$. Thus
1) follows from part a) of \maincurve 
 and from the fact that the images
of $\beta(3,1)$ and $\beta(2,2)$ are equal. 
2) follows likewise from part a)
 of \maincurve (note that in this case $\beta(2,2)$ surjects). Finally 
3) follows from part b) of Proposition.   
$\square$

\bigskip

\noindent (\mainpr \ns) {\it Proof of \main \ns:} The proof rests on
\maincurve \ns. The idea is ``to lift" the generators of
$R$ to the canonical ring of $S$. Let us define 
$$H^0(K_S^{\otimes s}) \otimes H^0(K_S^{\otimes r}) @>
\alpha(s,t) >> H^0(K_S^{\otimes s+t}) \ ,$$ and let also denote
$\alpha(s,1)$ as $\alpha_s$. As in the case of $R$, the images of
$\alpha(s,t)$ will tell us the generators of each graded piece
of the canonical ring of $S$. In fact it will suffice to prove the
following:

\medskip

\item{(a)} If $r=1$ and $n=2$, $\alpha_l$ surjects for all $l\geq
1$, except if $l =3$. The images of $\alpha_3=\alpha(3,1)$ and
$\alpha(2,2)$ are equal and have codimension $1$ in
$H^0(K_S^{\otimes 4})$.

\item{(b)} If $r=1$ and $n >2$, $\alpha_l$ surjects for all $l\geq
1$, except if $l =1,3$. The image of $\alpha_1$ has codimension
$n-2$ in $H^0(K_S^{\otimes 2})$. The map $\alpha(2,2)$ is
surjective. 

\item{(c)} If $r \geq 2$, $\alpha_l$ is surjective if
$l
\geq 3$. The image of $\alpha_1$ has codimension $r(n-2)$ in
$H^0(K_S^{\otimes 2})$. The image of $\alpha_2$ has codimension
$r-1$ in $H^0(K_S^{\otimes 3})$.

\medskip

Thus we proceed to prove (a), (b), (c). Recall that $Y$ is an
irreducible variety of minimal degree, and in particular, normal.
On the other hand the locus formed by the  points of $Y$ with non
finite fibers  has codimension $2$. Thus using  Bertini's
Theorem we can choose
 a smooth curve $C$ of $|K_S|$ such that the restriction of the
canonical morphism of $S$ to $C$ is finite (and flat) onto a smooth
rational normal curve of degree $r$. Let us denote by
$\theta$ the restriction of
$K_S$ to $C$. By adjunction $K_C =\theta^{\otimes 2}$. Since $K_S$
is base-point-free so is $\theta$. Finally, since $H^1(\Cal
O_X)=0$, $\pi$ is induced by the complete linear series $|\theta|$ and
therefore $C$, 
$\theta$ and $\pi$ satisfies the hypothesis of \maincurve \ns.

We prove first the statements in (a), (b) and (c) regarding the
maps $\alpha_l$. Consider the following commutative diagram:

$$\matrix
H^0((K_S^{\otimes l}) \otimes  H^0(\Cal O_S) &
\hookrightarrow & H^0((K_S^{\otimes l}) \otimes 
H^0(K_S)&
\twoheadrightarrow & H^0((K_S^{\otimes l}) \otimes 
H^0(\theta) \cr @VV   V @VV \alpha_l V
@VV  V\cr
 H^0(K_S^{\otimes l}) & \hookrightarrow &
H^0(K_S^{\otimes l+1})& \twoheadrightarrow &
H^0(\theta^{\otimes l+1})
\cr 
\endmatrix
$$

The right most horizontal arrows are
surjective because $H^1(\Cal O_S)=0$, by Serre duality and by Kawamata-Viehweg
vanishing.  The left hand side vertical arrow trivially surjects.
The right hand side vertical arrow is the composition of the map
$H^0(K_S^{\otimes l}) \otimes 
H^0(\theta) \longrightarrow H^0(\theta^{\otimes l}) \otimes 
H^0(\theta)$, which is surjective for all $l \geq 1$ again 
because $H^1(\Cal O_S)=0$, by Serre duality and by Kawamata-Viehweg vanishing,
and the map $\beta_l$ of multiplication of
global sections on
$C$,
studied in \maincurve \ns. Then it
follows from chasing the diagram that the map $H^0(K_S^{\otimes
l+1})
\longrightarrow  H^0(\theta^{\otimes l+1})$ maps the image of
$\alpha_l$ onto the image of $\beta_l$ and that the codimension of
the image of
$\beta_l$ in $H^0(\theta^{\otimes l+1})$ is equal to the
codimension of $\alpha_l$ in $H^0(K_S^{\otimes l+1})$. This,
together with \maincurve \ns, a.1, a.2, a.3, a.5 and b, proves the claims in
(a), (b) and (c) concerning the codimensions of the images of the
maps $\alpha_l$. 

Thus the only things left to prove are the claims about
$\alpha(2,2)$ when $r=1$. We consider now this commutative
diagram

$$\matrix
H^0(K_S^{\otimes 2}) \otimes  H^0(K_S) &
\hookrightarrow & H^0(K_S^{\otimes 2}) \otimes 
H^0(K_S^{\otimes 2})&
\twoheadrightarrow & H^0(K_S^{\otimes 2}) \otimes 
H^0(\theta^{\otimes 2}) \cr @VV \alpha_2  V @VV \alpha(2,2) V
@VV  V\cr
 H^0(K_S^{\otimes 3}) & \hookrightarrow &
H^0(K_S^{\otimes 4})& \twoheadrightarrow &
H^0(\theta^{\otimes 4})
\cr 
\endmatrix
$$ 
  
The right most horizontal arrows  are
surjective because $H^1(\Cal O_S)=0$ and by Serre duality, and by
Kawamata-Viehweg 
vanishing. The left hand side vertical arrow surjects, as we have
already proven. The right hand side vertical arrow is the
composition of the map
$H^0(K_S^{\otimes 2}) \otimes 
H^0(\theta^{\otimes 2}) \longrightarrow H^0(\theta^{\otimes 2}) \otimes 
H^0(\theta^{\otimes 2})$, which is surjective because $S$ is regular
and by Serre duality, and the map $\beta(2,2)$ of multiplication of
global sections on
$C$. Then it follows from chasing the diagram that
the map $H^0(K_S^{\otimes 4}) \twoheadrightarrow 
H^0(\theta^{\otimes 4})$ maps the image of $\alpha(2,2)$ onto the
image of $\beta(2,2)$ and that the codimension of the image of
$\beta(2,2)$ in $H^0(\theta^{\otimes 4})$ is equal to the
codimension of the image of $\alpha(2,2)$ in $H^0(K_S^{\otimes
4})$. On the other hand, we know that the image of $\beta(2,2)$ and of
$\beta_3=\beta(3,1)$ are equal of codimension $1$ in
$H^0(\theta^{\otimes 4})$, if $r=1$ and $n=2$. Thus  we conclude that 
the images of $\alpha(3,1)$ and $\alpha(2,2)$ in
$H^0(K_S^{\otimes 4})$ are also equal  and of codimension $1$. 
Finally, if $r=1$ and $n >2$, $\beta(2,2)$ surjects by \maincurve
\ns.a.4. Thus  we 
conclude that if $r=1$ and $n >2$, then $\alpha(2,2)$
surjects.  $\square$

\bigskip

\main complements 
 known results on
generation of the canonical ring of smooth, regular surfaces of
general type.  Ciliberto and Green (cf. \G \ns, Theorem 3.9.3, and
\Ci
\ns) proved that, given a
smooth surface of general type with $h^1(\Cal O_S)=0$ and
$K_S$ globally generated and being 
$\varphi$  the canonical morphism,   a
sufficient condition for the canonical ring of
$S$ to be generated in degree less than or equal to $2$ is that none
of the following happen:
\roster
\item $\varphi$ maps $S$ generically $2:1$ onto $\bold P^2$.
\item $\varphi(S)$ is a surface of minimal degree (other than $\bold
P^2$).
\endroster

As a corollary of Ciliberto and Green result and of \main 
we obtain the following

\proclaim{\maincor} Let $S$ be a smooth  regular surface of general
type and such that $K_S$ is globally
generated. Let
$\varphi$ be the canonical morphism of $S$. The canonical
ring of
$S$ is generated in degree less than or equal to $2$ if and only if
none of the following happens:
\roster
\item $\varphi$ maps $S$ generically $2:1$ onto $\bold P^2$.
\item $\varphi(S)$ is a surface of minimal degree (other than $\bold
P^2$).
\endroster
\endproclaim

\heading 3. Examples of surfaces of general type \endheading

In this section we construct some new examples of surfaces of general type
which satisfy the hypothesis of \main \ns. 
The easiest way one could think of producing examples 
would be to build suitable cyclic
covers of surfaces of minimal degree. However, as next remark shows,
only low degree cyclic covers can be induced by the canonical morphism
of a regular surface, so we have to employ other means to construct some
new examples.

\proclaim{\nocyc} Let $X$ be a surface 
of general type with at worst canonical singularities and with
base-point-free canonical bundle. 
Assume that the complete canonical series of $X$ restricts to a
complete linear series on a general hyperplane section (e.g., if $X$
is regular). Let $\varphi: X
\longrightarrow Y$ be the canonical morphism to a surface of minimal
degree. Let $n$ be the degree of $\varphi$. Assume that, on the complement
$U$ of a 
codimension $2$ closed subset of $Y$, 
$$f_*\Cal O_X = \Cal O_Y \oplus L^{-1} \oplus \cdots \oplus L^{\otimes
  1-n} \
.$$
Then $n=2$ or $3$. 
\endproclaim

{\it Proof}. Let $H$ be a general hyperplane section of $Y$ contained in
$U$ and let $C$ be the inverse image of $H$ by $\varphi$. Then $C$ is a smooth
irreducible member of $|K_X|$. By assumption the
morphism $\varphi|_C: C \longrightarrow H$ is induced by the complete linear
series of a line bundle $\theta$. By adjunction $\theta^{\otimes 2}
=K_C$. Thus $C$, $\theta$ and $\varphi|_C$ satisfy the hypothesis of \split
and $$(\varphi|_C)_*\Cal O_C=\Cal O_H \oplus (n-2)N^{-1}  \oplus N^{-2} \ 
.$$ 
On the other hand $(\varphi|_C)_*\Cal O_C$ is equal to the restriction of
$\varphi_*\Cal O_X$ to $H$, i.e, to 
$$ \Cal O_H \oplus (L \otimes \Cal O_H)^{-1} \oplus \cdots \oplus (L
\otimes O_H) ^{1-n} \
.$$ This is only possible if $n=2$ or $3$.
$\qed$

\proclaim{\nocyccor} Let $X$ be a regular surface of general type
with at worst canonical singularities and with base-point-free
canonical bundle. Let $Y$ be the image of $X$ by its canonical
morphism $X @> \varphi >> Y$. If $Y$ is a surface of minimal degree
and $\varphi$ 
is a cyclic cover, then the degree of $\varphi$ is $2$ or $3$. 
\endproclaim

The next proposition also rules out many possible examples of covers
of odd degree: 

\proclaim{\noodd} Let $X$ be a surface of general type with at
worst canonical singularities whose
canonical divisor is base-point-free. Let $\varphi$ be a morphism induced
by a subseries of $|K_X|$. If $\varphi$
is generically finite onto a smooth scroll $Y \subset \bold P^N$,
then the degree of $\varphi$ is even. In particular, there are not
generically finite covers of odd degree of smooth rational normal scrolls,
induced by subseries of $K_X$.  
\endproclaim

{\it Proof.}  
Let $f$ be a fiber of $Y$ and let $C$ be a section of $Y$. Let
$-d=C^2$. Since $Y$ is a scroll,
its hyperplane section is linearly equivalent to $C + mf$, for some
integer $m$. Then $K_X=\varphi^*(C+mf)$. Then $\text{ deg
  }\varphi=(\varphi^*f)\cdot (\varphi^*C)=(\varphi^*f)\cdot(K_X+\varphi^*f) $, which is an
even number. 
$\square$

\medskip

Now we mention some examples of regular minimal surfaces
$X$ whose
canonical morphism $\varphi$ maps onto a variety of
minimal degree and produce some new ones. 

\medskip

The cases when $\varphi$ is generically finite and has degree $2$
and $3$ have been completely studied by Horikawa and Konno (see \Hoone \ns,
Theorem 1.6, 
\Hotwo \ns, Theorem 2.3.I, 
\Hothree \ns, Theorem 4.1 and \Kon \ns, Lemma
2.2 and Theorem 2.3). As it turns out there exist generically double
covers of linear $\bold P^2$, the Veronese surface, smooth rational
normal scrolls $S(a,b)$ with $b \leq 4$ and cones over rational normal
curves of degree $2$, $3$ and $4$ and generically triple covers of
$\bold P^2$ (in particular cyclic triple covers of $\bold P^2$
ramified along a sextic with suitable singularities) and of the cones
over rational normal curves of degree $2$ 
and $3$.  Horikawa (see \Hofour \ns, 
Theorem 2.1) also describes all
generically finite quadruple 
covers $X @> \varphi >> 
Y$,  where $X$ is smooth, minimal regular surface, $\varphi$ is the
canonical morphism of $X$ and $Y$ is linear $\bold P^2$. 

\medskip

The examples of
of Horikawa and Konno just reviewed are examples of covers of degree less
than or equal to $3$ of surfaces of minimal degree and quadruple covers of
$\bold P^2$.
We now construct three new sets of examples of regular surfaces of
general type which are quadruple covers of surfaces
of minimal degree under the canonical morphism. These examples are 4:1
covers of smooth rational normal scrolls isomorphic to the Hirzebruch
surfaces $\bold F_0$ and $\bold F_1$ and of quadric cones in $\bold P^3$.

\medskip

\proclaim{\Fzero}
We construct  finite quadruple covers $X @> \varphi >> Y$, where $X$ is a
smooth minimal regular surface of general type, $\varphi$ is the canonical
morphism of $X$ and $Y$ is a smooth rational scroll $S(m,m)$, $m \geq
1$.
\endproclaim

 Let
$f$ be a fiber of one of the fibrations of $\bold P^1$ and let $f'$
be a fiber  of the other fibration. Then
$Y$ is  $\bold P^1 \times \bold P^1$ and it is embedded in $\bold
P^{2m+1}$ by $|f+mf'|$ or by $|f'+mf|$. 
If $Y$ is embedded by  $|f+mf'|$, let
$a_1,a_2,b_1$ and $b_2$ satisfy 
the following: either $a_1=1,  a_2=2, b_1=m+1$ and
$b_2=1$ or  
$a_1=2, a_2=1, b_1=1$ and $b_2=m+1$.
  If $Y$ is embedded by  $|f'+mf|$, let
$a_1,a_2,b_1$ and $b_2$ satisfy 
the following: either $b_1=1,  b_2=2, a_1=m+1$ and
$a_2=1$ or  
$b_1=2, b_2=1, a_1=1$ and $a_2=m+1$.
Let $D_i$s be smooth divisors linearly equivalent to $2(a_if+b_if')$
intersecting at $D_1 \cdot D_2$ distinct points. Those
divisors exist because by the choices of $a_1,a_2,b_1$ and $b_2$, both
$2(a_1f+b_1f')$ and $2(a_2f+b_2f')$ are base-point-free. 
Let $X' @> \varphi_1 >> Y$ be the
double cover of $Y$ ramified along $D_1$. Since $D_1$ is smooth, so is
$X'$. Let $D_2'$ be the inverse
image in $X'$ of $D_2$ by $\varphi_1$. Since $D_2$ is smooth and meets $D_1$
at distinct points, $D_2'$ is also smooth. Let $X @> \varphi_2 >> X'$
be the double cover of $X'$ ramified along $D_2'$. Since $X'$ and
$D_2'$ are both smooth, so is $X$. Let us call $\varphi=\varphi_1 \circ
\varphi_2$. Now we will show that $X$ is a regular surface
of general type, that $K_X = \varphi^* \Cal O_Y(1)$ and that $\varphi$ is
induced by the complete canonical series of $X$.
First we find out the structure of $\varphi_*\Cal O_X$ as module over
$\Cal O_Y$. Recall that ${\varphi_2}_*\Cal O_X=\Cal O_{X'} \oplus
{\varphi_1}^*\Cal O_Y(-a_2f-b_2f')$. Then $$\varphi_*\Cal O_X = {\varphi_1}_*\Cal
O_{X'} \oplus {\varphi_1}_*({\varphi_1}^*\Cal O_Y(-a_2f-b_2f'))\ .$$
Since ${\varphi_1}_*\Cal O_{X'} = \Cal O_Y \oplus \Cal O_Y(-a_1f-b_1f')$,
then by projection formula we have 
$$\varphi_*\Cal O_X = \Cal O_Y \oplus \Cal O_Y(-a_1f-b_1f') \oplus \Cal
O_Y(-a_2f-b_2f') \oplus \Cal O_Y(-(a_1+a_2)f-(b_1+b_2)f')  
\ .  $$
We see now that $X$ is regular. Recall that $H^1(\Cal O_X) =
H^1(\varphi_*\Cal O_X)$.  
Our choice of $a_1,a_2,b_1$ and $b_2$ implies that $a_1f+b_1f'$ and
$a_2f+b_2f'$ 
are both base-point-free and big, thus by Kawamata-Viehweg vanishing, 
$$\displaylines{H^1(\Cal O_Y(-a_1f-b_1f')) = H^1(\Cal
O_Y(-a_2f-b_2f'))= \cr H^1(\Cal O_Y(-(a_1+a_2)f-(b_1+b_2)f')) = 0 \ .}$$
Then, since $H^1(\Cal O_Y)$ also vanishes so does $H^1(\varphi_*\Cal O_X)$
and $H^1(\Cal O_X)$. 
We now compute $K_X$. Since $\varphi_2$ is a double cover ramified along
$D_2'$, $K_X=
{\varphi_2}^*(K_{X'} \otimes \varphi_1^*(\Cal O_Y(a_2f+b_2f'))$. By  a similar
reason, 
$K_{X'}=\varphi_1^*(K_Y \otimes \Cal O_Y(a_1f+b_1f'))$. Then 
$K_X = \varphi^*(K_Y \otimes \Cal O_Y((a_1+a_2)f+(b_1+b_2)f'))$. Since
$K_Y=\Cal O_Y(-2f-2f')$, it follows again from the choices of $a_1,
a_2, b_1$ and $b_2$ that
$K_X=\varphi^*\Cal O_Y(1)$. 
Finally, to see that $\varphi$ is induced by the
complete canonical linear series of $X$ we compute $H^0(K_X)$. We do
the computation in the case $\Cal O_Y(1) = \Cal O_Y(f+mf')$. The case
$\Cal O_Y(1) = \Cal O_Y(mf+f')$ is analogous. Since
$K_X=\varphi^*\Cal O_Y(1)$, 
$$\displaylines{H^0(K_X)= H^0(\Cal O_Y(1)) \oplus H^0(\Cal
  O_Y((1-a_1)f+(m-b_1)f')) 
\oplus  \cr H^0(\Cal 
O_Y((1-a_2)f+(m-b_2)f')) \oplus H^0(\Cal
O_Y((1-a_1-a_2)f+(m-b_1-b_2)f')) \ .}$$ Again, by the choices of $a_1,
a_2, b_1$ and $b_2$, the last three direct sums of the above
expression are $0$, so $\varphi$ is indeed induced by the complete
canonical series of $X$. 
$\qed$

\proclaim{\Fone} 
We construct finite quadruple covers $X @> \varphi >> Y$, where $X$ is a
smooth regular surface of general type with base-point-free canonical
bundle, $\varphi$ is the canonical 
morphism of $X$ and $Y$ is a smooth rational scroll $S(m-1,m)$, 
$m \geq 2$.
\endproclaim

Let
$C_0$ be a minimal section of $\bold F_1$ and let $f$ be one of the
fibers. 
Then 
$Y$ is  $\bold F_1$ and it is embedded in $\bold
P^{2m}$ by $|C_0+mf|$. 
Let $a_1,a_2,b_1$ and $b_2$ satisfy
the following:
either $a_1=1, a_2=2, b_1=m+1$ and $b_2=2$ or 
$a_1=2, a_2=1, b_1=2$ and $b_2=m+1$. 

Let $D_i$s be smooth divisors linearly equivalent to $2(a_iC_0+b_if)$
intersecting at $D_1 \cdot D_2$ distinct points. The fact that such
divisors exist
follows from our choice of $a_1, a_2, b_1$ and $b_2$, which implies that
the linear systems 
of $D_1$ and $D_2$ are base-point-free. 
Let $X' @> \varphi_1 >> Y$ be the
double cover of $Y$ ramified along $D_1$. Since $D_1$ is smooth, so is
$X'$. Let $D_2'$ be the inverse
image in $X'$ of $D_2$ by $\varphi_1$. Since $D_2$ is smooth and meets $D_1$
transversally, $D_2'$ is also smooth. Let $X @> \varphi_2 >> X'$
the double cover of $X'$ ramified along $D_2'$. Since $X'$ and
$D_2'$ are both smooth, so is $X$. Let us call $\varphi=\varphi_1 \circ
\varphi_2$. Now we will show that $X$ is a regular surface
of general type, that $K_X = \varphi^* \Cal O_Y(1)$ and that $\varphi$ is
induced by the complete canonical series of $X$.
First we find out the structure of $\varphi_*\Cal O_X$ as module over
$\Cal O_Y$. Recall that ${\varphi_2}_*\Cal O_X=\Cal O_{X'} \oplus
{\varphi_1}^*\Cal O_Y(-a_2C_0-b_2f)$. Then $$\varphi_*\Cal O_X = {\varphi_1}_*\Cal
O_{X'} \oplus {\varphi_1}_*({\varphi_1}^*\Cal O_Y(-a_2C_0-b_2f))\ .$$
Since ${\varphi_1}_*\Cal O_{X'} = \Cal O_Y \oplus \Cal O_Y(-a_1C_0-b_1f)$,
then by projection formula we have 
$$\displaylines{\varphi_*\Cal O_X = \Cal O_Y \oplus \Cal O_Y(-a_1C_0-b_1f) \oplus \Cal
O_Y(-a_2C_0-b_2f) \cr \oplus \Cal O_Y(-(a_1+a_2)C_0-(b_1+b_2)f)  
\ .}  $$
We see now that $X$ is regular. Recall that $H^1(\Cal O_X) =
H^1(\varphi_*\Cal O_X)$.  
Our choices of $a_1,a_2, b_1$ and $b_2$ imply that $a_1C_0+b_1f$ and
$a_2C_0+b_2f$ 
are both base-point-free and big divisors, thus by Kawamata-Viehweg vanishing, 
$$\displaylines{H^1(\Cal O_Y(-a_1C_0-b_1f)) = H^1(\Cal
O_Y(-a_2C_0-b_2f))= \cr H^1(\Cal O_Y(-(a_1+a_2)C_0-(b_1+b_2)f)) = 0 \ .}$$
Then, since $H^1(\Cal O_Y)$ also vanishes so does $H^1(\varphi_*\Cal O_X)$
and therefore $H^1(\Cal O_X)$. 
We now compute $K_X$. Since $\varphi_2$ is a double cover ramified along
$D_2'$, $K_X=
{\varphi_2}^*(K_{X'} \otimes \varphi_1^*(\Cal O_Y(a_2C_0+b_2f))$. By similar
reason, 
$K_{X'}=\varphi_1^*(K_Y \otimes \Cal O_Y(a_1C_0+b_1f))$. Then 
$K_X = \varphi^*(K_Y \otimes \Cal O_Y((a_1+a_2)C_0+(b_1+b_2)f))$. Since
$K_Y=\Cal O_Y(-2C_0 -3f)$, it follows from our choice of $a_1, a_2,
b_1$ and $b_2$ that
$K_X=\varphi^*\Cal O_Y(1)$. 
Finally, to see that $\varphi$ is induced by the
complete canonical linear series of $X$ we compute $H^0(K_X)$. Since
$K_X=\varphi^*\Cal O_Y(1)$, 
$$\displaylines {H^0(K_X)= H^0(\Cal O_Y(1)) \oplus H^0(\Cal O_Y((1-a_1)C_0+(m-b_1)f))
\oplus \cr H^0(\Cal 
O_Y((1-a_2)C_0+(m-b_2)f)) \oplus H^0(\Cal
O_Y((1-a_1-a_2)C_0+(m-b_1-b_2)f)) \ .}$$ Again, by the choices of $a_1,
a_2, b_1$ and $b_2$, the last three direct sums of the above
expression are $0$, so $\varphi$ is indeed induced by the complete
canonical series of $K_X$. 
$\qed$

\proclaim{Remark 3.6} With the same arguments, if we allow certain mild
singularities in $D_1$ and $D_2'$, one can construct
examples of covers of $\bold F_0$ and $\bold F_1$ with at worst
canonical singularities.
\endproclaim

Finally we construct an example of a quadruple cover of a singular
surface of minimal degree.

\proclaim{\cone} 
We construct an example of a smooth, generically finite,
 quadruple cover $X @> \varphi >> Z$ 
 of the quadric cone $Z$
in $\bold P^3$, where $X$ is a regular surface of general type whose
 canonical divisor is base-point-free, and $\varphi$ is its canonical
 morphism. 
\endproclaim

Let $Y=\bold F_2$. Let $C_0$ be the minimal section of $Y$ and let $f$
be a fiber of $Y$. Let $D_1$ be a smooth divisor on $Y$, lineraly
equivalent to $2C_0+6f$ and meeting $C_0$ transversally. Let $D_2$ be
a smooth divisor on $Y$ linerarly equivalent to $3C_0+6f$ and meeting
$D_1$ transversally. Such divisors
$D_1$ and $D_2$ exist because $2C_0+6f$ is very ample and $3C_0+6f$ is 
base-point-free. Note also that, since $(3C_0+6f) \cdot C_0=0$, $C_0$
and $D_2$ do not meet. Let $X' @> \varphi_1 >> Y$ be the double cover of
$Y$ along $D_1$. Since $D_1$ is smooth, so is $X'$. Since $D_1$ meets
$C_0$ at two distinct points, the 
pullback $C_0'$ of $C_0$ by $\varphi_1$ is a smooth line with
self-intersection $-4$. Let $D_2'$ be the pullback of $D_2$ by
$\varphi_1$. Since $D_1$ and $D_2$ meet transversally, $D_2'$ is smooth
and since $D_2$ and $C_0$ do not meet, neither do $D_2'$ and
$C_0'$. Let $L_2'$ be the pullback of $2C_0+3f$ by $\varphi_1$. 
Let $X @> \varphi_2 >> X'$ be the double cover of $X'$ along $D_2'
\cup C_0'$. Since $D_2'
\cup C_0'$ is smooth, so is $X'$. Let us denote $\varphi=\varphi_1
\circ \varphi_2$.  Then 
$$\displaylines{ \varphi_* \Cal O_X = {\varphi_1}_*{\varphi_2}_*\Cal
  O_X = {\varphi_1}_*(\Cal O_{X'} 
\oplus {\Cal L_2'}^*)= \cr \Cal O_Y \oplus \Cal O_Y(-C_0-3f) \oplus \Cal
O_Y(-2C_0-3f) \oplus \Cal O_Y(-3C_0-6f) \ . (3.7.1)}$$
Since $-C_0-3f$ and $-3C_0-6f$ are big and base-point-free, by
Kawamata-Viehweg vanishing and Serre duality, $H^1(\Cal
O_Y(-C_0-3f))=H^1(\Cal O_Y(-3C_0-6f))=0$. By Serre duality $H^1(\Cal
O_Y(-2C_0-3f))=H^1(\Cal O_Y(-f))^*=0$. Then, since $H^1(\Cal O_Y)=0$,
$X$ is regular.  
Arguing as in \Fzero and \Fone we see that
$$K_X=\varphi^*(K_Y \otimes \Cal O_Y(3C_0+6f))= \varphi^*\Cal O_Y(C_0+2f) (3.7.2) \ .$$
Now we compute $H^0(K_X)$. Using projection formula and 3.7.1 and
3.7.2 we
obtain that
$$\displaylines{H^0(K_X)=H^0(\Cal O_Y(C_0+2f)) \oplus H^0(\Cal O_Y(-f)) \oplus
H^0(\Cal O_Y(-C_0-f)) \oplus \cr
H^0(\Cal O_Y(-2C_0-4f)) = H^0(\Cal
O_Y(C_0+2f)) \ .}$$
Thus the canonical morphism of $X$ is the composition of $\varphi$ and the
morphism $Y @> \phi >> Z  \subset \bold P^3$, induced by the complete
linear series of $C_0+2f$. Since $\phi$ contracts $C_0$, the canonical
morphism of $X$ is not finite, but it is generically finite of degree
$4$ onto $Z$, which is a surface of minimal degree as we wanted. 

On the other hand, if $C_0''$ is the pullback of $C_0$ by $\varphi$, then
$C_0''$ is a smooth line with self-intersection $-2$. Thus the
morphism $\phi
\circ \varphi$ also factorizes as $\varphi' \circ \psi$, where
$$X @> \psi >> \overline X$$ is the morphism from $X$ to its canonical
model $\overline X$ and $$\overline X @> \varphi' >> Z$$ is the canonical
morphism of $\overline X$. Thus $\varphi'$ is an example of a finite,
$4:1$ canonical morphism from regular surface of general type with
canonical singularities onto a singular surface of minimal
degree. $\qed$

\heading 4. The canonical ring of a curve \endheading

In this section we study a very well known ring, the
canonical ring of a curve. It was proved by Noether that the canonical
ring of smooth curve $C$ is generated in degree $1$ if and only if the curve is
non hyperelliptic. Our purpose is to show again the link between the
structure of this canonical ring and the structure of $\Cal O_C$ as an algebra
over $\Cal O_{\bold P^1}$ via a suitable morphism from $C$ to $\bold
P^1$. 
Precisely, we will look at curves endowed with certain finite
morphisms to a rational curve. First, we will consider 
hyperelliptic curves and its $2:1$ morphism to $\bold
P^1$. Second, we will consider curves having a base-point-free
theta-characteristic inducing a morphism from the curve to a rational
normal curve. The latter class of curves include hyperelliptic curves
but also other curves far more general: those curves with a
base-point-free theta characteristic with two linearly independent sections.  
We  will give a new proof of Noether's theorem in these two cases. 

\proclaim{\hyp} Let $C$ be a smooth hyperelliptic
curve of genus $g$. 
\item{1)} If $g=2$, then its canonical ring is generated by its
  elements of degree
$1$, and by  
$1$ element of degree $3$. 
\item{2)} If $g \geq 3$, then its canonical ring is generated by its
  elements of degree
$1$, and by  
$g-2$ elements of degree $2$.
\endproclaim

{\it Proof}. 
Let $L$ be the base-point-free $g^1_2$ on $C$ and let $\pi: C
\longrightarrow \bold P^1$ be the morphism induced by
$|L|$. The cover $\pi$ is double, therefore 
$\pi_*\Cal O_C = \Cal O_{\bold P^1} \oplus \Cal O_{\bold
P^1}(-g-1)$ and the structure of $\Cal O_{\bold P^1} \oplus \Cal O_{\bold
P^1}(-g-1)$  as algebra over $\Cal O_{\bold P^1}$ is well known: 
the elements of $\Cal O_{\bold P^1}$ multiply among themselves by the
product in $\Cal O_{\bold P^1}$, the elements of $\Cal O_{\bold P^1}$
multiply the elements of $\Cal O_{\bold
P^1}(-g-1)$ via the module structure of $\Cal O_{\bold
P^1}(-g-1)$ over $\Cal O_{\bold P^1}$, and finally the multiplication
of elements of
$\Cal O_{\bold 
P^1}(-g-1)$ is dictated by the branch divisor on $\bold P^1$ and
produces elements of $\Cal O_{\bold P^1}$.     Recall now that
$K_C=L^{\otimes g-1}$.
Then $\pi_*K_C^{\otimes n}= \Cal O_{\bold P^1}(n(g-1)) \oplus \Cal O_{\bold
  P^1}((n-1)(g-1)-2)$. This implies for instance that $H^0(K_C) =
H^0(\Cal O_{\bold 
  P^1}(g-1))$ and $H^0(K_C^{\otimes 2}) = H^0(\Cal O_{\bold
  P^1}(2(g-1))) \oplus H ^0(\Cal O_{\bold
  P^1}(g-3))$. Thus 
the image of $$H^0(K_C) \otimes H^0(K_C) @> \alpha_1 >>
H^0(K_C^{\otimes 2}) $$ is  $H^0(\Cal O_{\bold
  P^1}(2(g-1)))$, which has codimension $g-2$ in $H^0(K_C^{\otimes
  2})$. On the other hand, since $H^0(\Cal O_{\bold
  P^1}((n-1)(g-1)-2)) \neq 0$ if $n \geq 2$, except for $n=g=2$, we
see that the map 
$$H^0(K_C^{\otimes n}) \otimes H^0(K_C) @> \alpha_n >>
H^0(K_C^{\otimes n+1}) $$
surjects for all $n \geq 2$, except if $n=g=2$. 
If $g=2$, $\alpha_2$ does not surject, since its image is $H^0(\Cal
O_{\bold P^1}(3))$, which has codimension $1$ in $H^0(K_C^{\otimes
  3})$. $\qed$ 

\proclaim{\Noet} Let $C$ be a smooth curve of genus $g \geq 3$
possesing a base-point-free line bundle $\theta$ such that 
$\theta^{\otimes 2}=K_C$ and such that $|\theta|$ induces a morphism
of degree $n$
onto a rational normal curve of degree $r$. 
Then the canonical ring of $C$ is generated
in degree $1$ unless $r=\frac{g-1}{2}$, equivalently, unless $n=2$. In
particular, if $r \neq \frac{g-1}{2}$, then $C$ is  non-hyperelliptic. 
\endproclaim 
 
{\it Proof.}
Let $\pi$ be the morphism induced by $|\theta|$. 
The pair $(C,\theta)$ satisfies the hypothesis of \split \ns. Thus
we know that 
$$E= (n-2)\Cal O_{\bold P^1}(-r-1) \oplus \Cal O_{\bold
P^1}(-2r-2) \ .
$$
Recall also (cf. \algstr \ns) that
$$\displaylines{ \text{ For
some $1 \leq i,j \leq n-2$, the projection of the map } \cr
\Cal
O_{\bold P^1}(a_i) \otimes \Cal O_{\bold P^1}(a_j) \longrightarrow
\pi_*\Cal O_C \text{ to }\Cal O_{\bold P^1}(-2r-2) \text{ is
surjective,} \cr 
\text{ in fact, it is an isomorphism.} }$$   
We want to discuss the surjectivity of the maps
$$H^0(K_C^{\otimes m}) \otimes H^0(K_C) @> \gamma_m >> 
H^0(K_C^{\otimes m+1}) $$
for all $m \geq 1$. The map $\gamma_m$ is in fact the map
$\beta(2m,2)$ defined in \maincurve \ns. Recall that 
$$\displaylines{H^0(K_C)=H^0(\Cal O_{\bold P^1}(2r)) \oplus
  (n-2)H^0(\Cal O_{\bold P^1}(r-1)) \text{ and }\cr 
H^0(K_C^{\otimes 2})= H^0(\Cal
O_{\bold P^1}(4r))
\oplus (n-2)H^0(\Cal O_{\bold P^1}(3r-1)) \oplus H^0(\Cal
O_{\bold P^1}(2r-2))\ .}$$ 
Then the restriction of $\beta(2,2)$ 
to $H^0(\Cal O_{\bold P^1}(2r)) \otimes H^0(\Cal O_{\bold P^1}(2r))$
surjects onto \linebreak
$H^0(\Cal
O_{\bold P^1}(4r))$. Thus if $n=2$, the image of $\beta(2,2)$ is $H^0(\Cal
O_{\bold P^1}(4r))$, which has codimension $2r-1$. Note that, since in
this case $C$ is hyperelliptic, $r = \frac {g-1}{2}$ 
(and $g$ is therefore odd), so the codimension of the image of
 $\beta(2,2)$ in $H^0(K_C^{\otimes 2})$ is $g-1$, as seen in
 \hyp \ns. 
If $n \geq 3$, then restriction of $\beta(2,2)$ to 
 $H^0(\Cal O_{\bold P^1}(2r)) \otimes (n-2)H^0(\Cal O_{\bold
   P^1}(r-1))$ surjects onto 
$(n-2)H^0(\Cal O_{\bold P^1}(3r-1))$. Finally, it follows from \algstr 
that the image of the restriction of $\beta(2,2)$ to $(n-2)H^0(\Cal O_{\bold
   P^1}(r-1)) \otimes (n-2)H^0(\Cal O_{\bold
   P^1}(r-1))$ projects onto $H^0(\Cal O_{\bold
P^1}((2r-2))$. Therefore $\beta(2,2)$ surjects if $n \geq 3$. 
On the other hand
$$\displaylines{H^0(K_C^{\otimes m})= H^0(\Cal
O_{\bold P^1}(2mr))
\oplus \cr (n-2)H^0(\Cal O_{\bold P^1}((2m-1)r-1)) \oplus H^0(\Cal
O_{\bold P^1}((2m-2)r-2)) \ ,}$$
therefore $\beta(2m,2)$ surjects for all $m \geq 2$. $\qed$

\medskip 

Not every curve $C$ of genus $g$ has a theta-characteristic satisfying
the hypothesis of \Noet \ns. The curves with
theta-characteristics with a positive even number of sections form a
divisor $\Cal M_g^1$ in $\Cal M_g$ (see \Be and \F \ns). 
 Moreover, the
theta-characteristic of  a general curve
of $\Cal M_g^1$ is base-point-free (see \T \ns). Thus we can deduce
from this and the above proposition Noether's theorem for the case in
which $C$ is a general 
curve of $\Cal M_g^1$. 

\proclaim{\Noetcor} Let $C$ a smooth curve of genus $g \geq 3$,
general in $\Cal M_g^1$. Then  the canonical ring of $C$ is generated
in degree $1$. In particular, $C$ is non-hyperelliptic. 
\endproclaim

\heading 5. Homogeneous rings of Calabi-Yau threefolds \endheading

In this section  we give a new, different proof of the following
 result contained in \CY as
part of Theorems 1.4 and 1.7 and Corollary 1.8 of that article. The
arguments we use will give yet another illustration of how
the algebra structure induced on $\Cal O_X$ is linked to the very ampleness
and normal generation of a line bundle.

\proclaim{\CYth}
Let $X$ be  a Calabi-Yau threefold and let $B$ an ample
and base-point-free line bundle such that $h^0(B)=4$. Let $S$ be  a divisor in
$|B|$ and let $C$ be any smooth curve $C \in |B \otimes
\Cal O_S|$. 
The following are
 equivalent: 
\roster 
\item $B^{\otimes 2}$ satisfies property $N_0$.
\item
$B^{\otimes 3}$ satisfies property $N_0$. 
\item The sectional genus of $B$ is bigger than $3$.
\item The curve $C$ is non-hyperelliptic. 
\endroster

\endproclaim

{\it Proof.}
We will denote by $\varphi$  the morphism
induced by $|B|$ onto
$\bold P^3$. First of all we observe that, by \M \ns, Theorem 2,  
$$H^0(B^{\otimes n}) \otimes H^0(B) \longrightarrow H^0(B^{\otimes
  n+1}) $$
surjects if $n \geq 4$. This implies that
$$H^0(B^{\otimes n}) \otimes H^0(B^{\otimes n}) \longrightarrow
H^0(B^{\otimes 2n}) $$
surjects if $n \geq 4$. Therefore, $B^{\otimes 2}$ satisfies property
$N_0$ if and only if 
$$H^0(B^{\otimes 2}) \otimes H^0(B^{\otimes 2}) \longrightarrow
H^0(B^{\otimes 4}) $$ surjects. 
Analogously,  $B^{\otimes 3}$ satisfies property
$N_0$ if and only if 
$$H^0(B^{\otimes 3}) \otimes H^0(B^{\otimes 3}) \longrightarrow
H^0(B^{\otimes 6}) $$ surjects. 
The proof goes on through three steps.
\vskip .2 truecm
\noindent {\bf (5.1.1)} {\it The vector bundle $\varphi_*\Cal
O_X$}.  The morphism $\pi$ is finite of degree
$n \geq 2$.  The push down of $\Cal O_X$ by  $\varphi$,
$\varphi_*\Cal O_X$, is isomorphic to $\Cal O_{\bold P^3}
\oplus E$, where $E$ is a vector bundle of rank $n-1$ on $\bold P^3$.  
Since  $h^1(B^{\otimes n})=h^2(B^{\otimes n})=0$ for all
$n \in \bold Z$, using projection formula we obtain from
Horrocks' criterion that $E$ splits as a direct sum of line
bundles. On the other hand $B \otimes \Cal O_S=K_S$, and the
restriction of $|B|$ is the complete canonical series of $S$, for
$h^1(\Cal O_X)=0$. In addition, $h^1(\Cal O_S)=0$. 
Let us denote by $\pi$ the restriction of $\varphi$ to $C$. We are under
the hypothesis of \split \ns, therefore
$$\pi_*\Cal O_C = \Cal O_{\bold P^1} \oplus \Cal (n-2)O_{\bold
  P^1}(-2) \oplus \Cal O_{\bold P^1}(-4) \ .$$
This implies that $$\varphi_*\Cal
O_X= \Cal O_{\bold P^3} \oplus (n-2)\Cal O_{\bold
  P^3}(-2) \oplus \Cal O_{\bold P^3}(-4) \ .$$
Let us call $E_1=(n-2)\Cal O_{\bold
  P^3}(-2)$ and $E_2=\Cal O_{\bold P^3}(-4)$.

\vskip .2 truecm
\noindent {\bf (5.1.2)} {\it Relationship between the algebra structure
  of $\varphi_*\Cal 
O_X$ and the normal generation of $B^{\otimes 2}$ and
$B^{\otimes 3}$}.  We study now the
$\Cal O_{\bold P^3}$-algebra structure of $\varphi_*\Cal
O_X$ in relation with the surjectivity of the maps
$$H^0(B^{\otimes 2}) \otimes H^0(B^{\otimes 2}) @>\alpha >>
H^0(B^{\otimes 4})$$ and  
$$H^0(B^{\otimes 3}) \otimes
H^0(B^{\otimes 3}) @>\beta >> H^0(B^{\otimes 6})\ .$$ Recall that
$B^{\otimes 2}$ (resp. $B^{\otimes 3}$) is normally generated if and
only if $\alpha$ (resp. $\beta$) surjects. 
Recall that $H^0(B^{\otimes 2})=H^0(\Cal O_{\bold P^3}(2))
\oplus H^0(E(2))=H^0(\Cal O_{\bold P^3}(2)) \oplus
H^0(E_1(2))$ and $H^0(B^{\otimes 4})= H^0(\Cal O_{\bold
P^3}(4))
\oplus H^0(E(4))=H^0(\Cal O_{\bold
P^3}(4))
\oplus H^0(E_1(4))\oplus H^0(E_2(4))$. Then we can see
the map $\alpha$ as direct sum of:

$$\displaylines{ H^0(\Cal O_{\bold P^3}(2)) \otimes
H^0(\Cal O_{\bold P^3}(2)) @> \gamma >> H^0(\Cal O_{\bold
P^3}(4)) \cr
H^0(\Cal O_{\bold P^3}(2)) \otimes H^0(E_1(2)) @> \delta
>> H^0(E_1(4)) \cr
H^0(E_1(2)) \otimes H^0(\Cal O_{\bold P^3}(2))  @> \epsilon >>
H^0(E_1(4)) \cr
H^0(E_1(2)) \otimes H^0(E_1(2))  @> \eta >>H^0(\Cal O_{\bold P^3}(4)
\oplus E(4)) \ .}$$

The map $\gamma$ is induced by ring multiplication on $\Cal
O_{\bold P^3}$ and it is therefore surjective. The maps
$\delta$ and $\epsilon$ are induced by module
multiplication
and are also surjective. Therefore $\alpha$ surjects if and
only if the composition of $\eta$ with the projection to $H^0(E_2(4))$ 
is surjective. 
Now the map $\eta$ depends on the way in which  elements of
$E$ multiply among themselves. Let us denote by $\mu$ the  morphism 
$E
\otimes E$ to $\Cal O_{\bold P^3} \oplus E$ induced by the ring
structure of $\varphi_*\Cal O_X$. Now, the composition of $\eta$ with the
projection to $H^0(E_2(4))$  
is surjective if and only if 

\vskip .1 truecm
\noindent (*) $\mu$ induces
an isomorphism from at least one of the  components of $E_1
\otimes E_1$ isomorphic to $\Cal O_{\bold P^3}(-2) \otimes
\Cal O_{\bold P^3}(-2)$ onto the component $E_2$. 
\vskip .1 truecm

The same argument proves that $\beta$ is surjective  if
and only if (*) holds. 
\vskip .2 truecm
\noindent {\bf (5.1.3)} {\it Algebra structure of $\varphi_*\Cal
O_X$ and the curve $C$}. On the other hand, \linebreak
$(\varphi|_C)_*(\Cal O_C)$ and
$\Cal O_{\bold P^1} \oplus (E \otimes \Cal O_{\bold P^1})$ are
isomorphic as $\Cal  
O_{\bold P^1}$-algebras,  the former with the algebra
structure induced by the cover $\varphi|_C: C \longrightarrow \bold P^1$
of $\bold P^1$  and
the latter with the algebra structure inherited from the
algebra structure of $\varphi_*\Cal O_X$. Arguing similarly
as before, the multiplication map $H^0(K_C)
\otimes H^0(K_C) @>\lambda >> H^0(K_C^{\otimes 2})$ can be
related  to the algebra structure of $(\varphi|_C)_*(\Cal O_C)$, and in
fact, $\lambda$ is surjective if and only if (*) holds. The
surjectivity of $\lambda$ is equivalent to $C$ being
non-hyperelliptic. 

Finally, the failure of (*) to hold is equivalent to the
fact that
  $\mu$ restricted to $E_1 \otimes
E_1$ projects to $0$ in $\Cal O_{\bold
P^3}(-4)$. This condition implies that $\varphi_*\Cal O_X$
contains a subalgebra of rank $(n-1)$, namely $\Cal O_{\bold
P^3} \oplus E_1$ and
hence
$\varphi$ decomposes as a cover of degree bigger than or
equal to
$2$ and a cover of degree $n-1$, the latter integral, as $C$
is. Therefore this is only possible  if
$n=2$, which is equivalent to
$g(C)=3$. $\square$

\heading References \endheading

\item{\Be} A. Beauville, {\it Prym varieties and the Schottky
    problem}, Invent. Math. {\bf 41} (1977), 149--196.
 
\item{\Bo} E. Bombieri, {\it Canonical models of surfaces of general type},
Inst. Hautes. Et. Sci. Publ. Math. {\bf 42} (1973), 171--219.

\item{\Ca} F. Catanese,{\it  On the moduli spaces of surfaces
of general type}, J. Differential Geom. {\bf 19} (1984),
483--515.

\item{\Ci} C. Ciliberto, {\it Sul grado dei generatori dell'anello di
    una superficie di tipo generale}, 
Rend.
Sem. Mat. Univ. Politec. Torino {\bf 41} (1983) 
 
\item{\EH} D. Eisenbud and J. Harris, {\it On varieties of minimal
    degree (A centennial account)}, Algebraic Geometry, Bowdoin 1985,
  Amer. Math. Soc. Symp. in Pure and App. Math. {\bf 46} (1987), 1--14.
 
\item{\F}  H.M. Farkas, {\it Special divisors and analytic subloci of
Teichmueller space},  Amer. J. Math. {\bf 88} 1966 881--901.

\item{\sos} F.J. Gallego and  B.P. Purnaprajna, {\it Projective normality and
syzygies of algebraic surfaces}, J. Reine Angew. Math. {\bf 506} (1999),
145--180.

\item{\CY} F.J. Gallego and B.P. Purnaprajna, {\it Very ampleness and higher
syzygies for Calabi-Yau threefolds}, Math. Ann. {\bf 312} (1998),
no. 1, 133--149.

\item{\can} F.J. Gallego and B.P. Purnaprajna, {\it Canonical Covers
    of varieties 
of minimal degree}, in preparation. 

\item{\trig} F.J. Gallego and B.P. Purnaprajna, {\it On the rings of
    trigonal curves}, 
in preparation.

\item{\G} M.L. Green, {\it The canonical ring of a variety of
general type}, Duke Math. J. {\bf 49} (1982), 1087--1113.

\item{\HM} D. Hahn and R. Miranda, {\it Quadruple covers of algebraic
varieties}, J. Algebraic Geom. {\bf 8} (1999), 1--30.

\item{\Hoone} E. Horikawa, {\it Algebraic surfaces of general type with
  small $C^2_1.$ I}, Ann. of Math. (2) {\bf 104} (1976),
  357--387. 

\item{\Hotwo} E. Horikawa, {\it Algebraic surfaces of general type with small
$c\sp{2}\sb{1}$, II}, Invent. Math. {\bf 37} (1976),  121--155.

\item{\Hothree} E. Horikawa, {\it Algebraic surfaces of general type with small
$c\sp{2}\sb{1}$, III}, Invent. Math. {\bf 47} (1978),  209--248.

\item{\Hofour} E. Horikawa, {\it E. Algebraic surfaces of general
    type with small
$c\sp{2}\sb{1}$, IV}, Invent. Math. {\bf 50} (1978/79),  103--128.

\item{\Kod} K. Kodaira, {\it Pluricanonical systems on algebraic surfaces of
general type}, J. Math. Soc. Japan {\bf 20} (1968), 170--192.

\item{\Kon} K. Konno, {\it Algebraic surfaces of general type with $c\sp 2\sb
1=3p\sb g-6$}, Math. Ann. {\bf 290} (1991), 77--107.

\item{\M} D. Mumford, {\it Varieties defined by quadratic equations},
  Corso CIME in Questions on Algebraic Varieties, Rome (1970),
  30--100. 

\item{\OP} K. Oguiso and T. Peternell, {\it On polarized canonical Calabi-Yau
threefolds}, Math. Ann. {\bf 301} (1995), 237--248. 

\item{\T} M. Teixidor i Bigas, {\it Half-canonical series on
algebraic curves}, Trans. Amer. Math. Soc. {\bf 302} (1987),  99--115.

\end